\documentclass[pdflatex,sn-mathphys-num]{sn-jnl}

\usepackage{graphicx}%
\usepackage{multirow}%
\usepackage{amsmath,amssymb,amsfonts}%
\usepackage{amsthm}%
\usepackage{mathrsfs}%
\usepackage[title]{appendix}%
\usepackage{xcolor}%
\usepackage{textcomp}%
\usepackage{manyfoot}%
\usepackage{booktabs}%
\usepackage{algorithm}%
\usepackage{algorithmicx}%
\usepackage{algpseudocode}%
\usepackage{listings}%
\usepackage{tikz-cd}

\newcommand*{\R}{\mathbb R}

\newcommand*{\upto}{\upharpoonright}

\newcommand*{\Eqvl}{\Longleftrightarrow}

\newcommand*{\ot}{\mathrm{ot}}

\newcommand*{\Dim}{\mathrm{Dim}}
\newcommand*{\cantor}{2^\omega}
\newcommand*{\ext}{\mathrm{Ext}}
\newcommand*{\strings}{2^{<\omega}}

\newcommand*{\brk}[1]{\text{\textlbrackdbl}#1\text{\textrbrackdbl}}

\theoremstyle{thmstyleone}%
\newtheorem{theorem}{Theorem}[section]

\theoremstyle{thmstyletwo}%

\theoremstyle{thmstylethree}%
\newtheorem{definition}[theorem]{Definition}%
\newtheorem{lemma}[theorem]{Lemma}

\newtheorem{corollary}[theorem]{Corollary}
\newtheorem{principle}[theorem]{Principle}
\newtheorem{question}[theorem]{Question}
\newtheorem{proposition}[theorem]{Proposition}%

\begin{document}

\title[]{On the Hausdorff dimension of maximal chains and antichains of Turing and Hyperarithmetic degrees}


\author*[1]{\fnm{Sirun} \sur{Song}}\email{songsirun@foxmail.com}

\author[1]{\fnm{Liang} \sur{Yu}}


\affil*[1]{\orgdiv{Math school}, \orgname{Nanjing university}, \orgaddress{\street{Hankou street}, \city{Nanjing}, \postcode{210089}, \state{Jiangsu}, \country{China}}}




\abstract{This paper investigates the Hausdorff dimension properties of chains and antichains in Turing degrees and hyperarithmetic degrees. Our main contributions are threefold: First, for antichains in hyperarithmetic degrees, we prove that every maximal antichain necessarily attains Hausdorff dimension 1. Second, regarding chains in Turing degrees, we establish the existence of a maximal chain with Hausdorff dimension 0. Furthermore, under the assumption that $\omega_1=(\omega_1)^L$, we demonstrate the existence of such maximal chains with $\Pi^1_1$ complexity.  Third, we extend our investigation to maximal antichains of Turing degrees by analyzing both the packing dimension and effective Hausdorff dimension. 
}




\maketitle

\section{Introduction}\label{sec:introduction}


In recursion theory, a set of Turing degrees is defined as a chain if every two distinct elements are Turing comparable, and is considered maximal if it cannot be properly extended. Conversely, a set of Turing degrees is an antichain if any two distinct elements are Turing incomparable, with maximality defined analogously. These concepts extend naturally to hyperarithmetic degrees.

In 2006, Liang Yu \cite{yu2006measure} demonstrated the existence of a non-measurable antichain of Turing degrees. This led Jockusch to inquire whether every maximal antichain of Turing degrees is non-measurable. A significant advancement came in 2015 when C.T. Chong and Liang Yu \cite{chong2016measure} resolved this question by constructing maximal antichains of both Turing degrees and hyperarithmetic degrees with Lebesgue measure 0.

A natural extension proposed by Liang Yu concerns the minimal possible Hausdorff dimension of maximal antichains in Turing degrees. While it is straightforward to verify that the antichains constructed by Chong and Yu attain Hausdorff dimension 1, determining this dimension in general remains challenging. This paper systematically investigates the Hausdorff dimension of chains and antichains in both Turing degrees and hyperarithmetic degrees. 

The outline of this paper is as follows: section 2 provides preliminary materials and demonstrates that the Chong-Yu antichains achieve Hausdorff dimension 1. In Section 3, we analyze maximal antichains of hyperarithmetic degrees. Inspired by the proof of higher Demuth theorem in \cite{chong2015randomness}, we prove that every such antichain attains Hausdorff dimension 1. Our approach leverages the theory of higher randomness (specifically, $\Pi^1_1$-random reals) and the critical property that these reals are $\Delta^1_1$-dominated.
  
Section 4 focuses on maximal chains of Turing degrees. Building on the observation that recursively traceable degrees have effective dimension 0, we prove a key results: any ascending countable sequence of recursively traceable degrees has a recursively traceable minimal cover. This enables the inductive construction of a maximal $\omega_1$-chain consisting solely of recursively traceable degrees, thereby establishing its effective Hausdorff dimension 0. Furthermore, under the axiom $\omega_1=(\omega_1)^L$, we show the existence of $\Pi^1_1$-definable maximal chains. 

Section 5 examines maximal antichains of Turing degrees. Unlike the hyperarithmetic case, the absence of universal recursive domination in 1-random reals (the hyperimmune-free property) precludes analogous strong results. Nevertheless, we prove that every maximal antichain has effective Hausdorff dimension 1, a result non-relativizable except for K-trivial oracles. By transitioning to packing dimension, however, we achieve a complete characterization: all maximal antichains attain packing dimension 1. This directly follows from Downey and Greenberg's theorem \cite{downey2008turing} on the packing dimension of minimal degree reals.

\section{Preliminaries}\label{sec:preliminaries}

We assume the reader is familiar with the recursion theory, algorithmic randomness theory, and the hyperarithmetic theory. For those who are not, they can refer to \emph{Turing computability: theory and applications} by Soare \cite{soare2016turing}, \emph{Algorithmic randomness and complexity} by Downey and Hirschfeldt \cite{downey2010algorithmic}, and \emph{Recursion theory: Computational aspects of definability} by C.T.Chong and Liang Yu \cite{chong2015recursion}.

Throughout this paper, we use sets to stand for sets of integers, classes to mean subsets of $\cantor$. The former is often written in capital English letters $X,Y,Z,\cdots$, the latter is written in calligraphic letters $\mathcal A,\mathcal B,\cdots$. Especially, $\mathcal C=\{X\in\cantor: X\in L_{\omega_1^X}\}$ is the class of the quickly constructive reals, where $L$ is the hierarchy of G\"odel's constructiable universe. 
Without making any confusion, we also use capital letters to denote infinite binary sequence, i.e. elements in $\cantor$. The set of finite binary sequences is denoted $\strings$, and we use $\sigma,\tau,\cdots$ to denote finite strings. Particularly, we use $\lambda$ to denote the empty string. Given a set $X\in\cantor$ and a string $\sigma$, we use $\sigma\prec X$ to denote that $\sigma$ is a prefix of $X$, we say $X$ extends $\sigma$. For any string $\sigma$, define the class $\brk{\sigma}=\{X\in\cantor:\sigma\prec X\}$, and for a set of strings $D\subseteq\strings$, let $\brk{D}=\{X\in\cantor: \exists \sigma\in D(\sigma\prec X)\}$. We equip the product topology on the space $\cantor$, and define the Lebesgue measure on it by letting $\mu(\brk{\sigma})=2^{-|\sigma|}$ for every $\sigma\in \strings$. There are other definitions of measures on $\cantor$, they can be defined over a premeasure, i.e. a function $\rho:\strings\to \mathbb R^{\geq 0}$ such that $\rho(\lambda)=0$ and  $\rho(\sigma)=\rho(\sigma0)+\rho(\sigma1)$ for every $\sigma\in\strings$. If the premeasure is recursive, i.e. the set $$\{(\sigma,p,q):p,q\in\mathbb Q\wedge \sigma\in \strings\wedge p<\rho(\sigma)<q\}$$is recursive, then we say the induced measure $\nu$ is a computable measure.

\begin{definition}
	\label{def:chains}
	\begin{enumerate}
		\item A class $\mathcal A\subseteq \cantor$ is called a chain of Turing degrees if \begin{enumerate}
			\item $\forall X\in \mathcal A\forall Y\in\cantor[Y\equiv_T X\to Y\in \mathcal A]$;
			\item $\forall X,Y\in \mathcal A[X\not\equiv_T Y\to X<_T Y\vee Y<_T X]$.
		\end{enumerate}
		A chain $\mathcal A$ is maximal if for any $Z\notin \mathcal A$, the class $\mathcal A\cup\{X\in\cantor:X\equiv_T Z\}$ is not a chain any more.
		\item A class $\mathcal B\subseteq \cantor$ is called an antichain of Turing degrees if \begin{enumerate}
			\item $\emptyset\notin \mathcal B$;
			\item $\forall X\in \mathcal A\forall Y\in\cantor[Y\equiv_T X\to Y\in \mathcal A]$;
			\item $\forall X,Y\in \mathcal A[X\not\equiv_T Y\to X\not<_T Y\wedge Y\not<_T X]$.
		\end{enumerate}
		An antichain $\mathcal B$ is maximal if for any $Z\notin \mathcal B$, the class $\mathcal B\cup\{X\in\cantor:X\equiv_T Z\}$ is not an antichain any more.
	\end{enumerate}
\end{definition}

The subject of chains and antichains probably goes back to Sacks \cite{sacks1966degrees}. In 2006, Liang Yu \cite{yu2006measure} proved the following result that revealed the measure property of the antichains of Turing degrees.

\begin{theorem}[Liang Yu \cite{yu2006measure}]
	For every locally countable partial order $\mathbb P =(\cantor,\leq_P)$, there is a non-measurable antichain in $\mathbb P$. In particularly, there is a non-measurable antichain in $(\cantor,\leq_T)$.
\end{theorem} 

Jockush asked if every maximal antichain of Turing degree is non-measurable. In the paper of C.T.Chong and Yu, they proved that the answer is negative. They showed that there are null maximal antichains both in Turing degrees and hyperdegrees.

\begin{theorem}[C.T.Chong and L.Yu \cite{chong2016measure}]
	\label{thm:Yu's-maximal-antichain}
	\begin{enumerate}
		\item There is a null maximal antichain $\mathcal A$ of hyperdegrees. Moreover, for any $\Pi^1_1$-random real $R$, there is an $X\in \mathcal A$ such that $R<_h X$.
		\item There is a null maximal antichain $\mathcal B$ of Turing degrees. Moreover, for any $\Pi^1_1$-random real $R$, there is an $X\in \mathcal B$ such that $R<_T X$ and $X\not\leq_h R$.
	\end{enumerate}
\end{theorem}

A set $R\in\cantor$ is said to be $\Pi^1_1$-random if $X$ does not belong to any $\Pi^1_1$ class of Lebesgue measure 0. Similarly, $R$ is $\Delta^1_1$-random if it doesn't belong to any $\Delta^1_1$ class of measure 0. The notion of $\Delta^1_1$-random can be viewed as a counterpart of the Schnorr random, by the following result.

\begin{proposition}[See \cite{chong2015recursion}]
	\label{prop:delta11-random-vs-schnorr}
	Suppose that $\mathcal A\subset\cantor$ is a $\Delta^1_1$-null class, then there exists a $\Delta^1_1$ set $V\subseteq \omega\times\strings$ such that $\forall n(\mu( \brk{V_n})=2^{-n})$, and $\mathcal A\subseteq \bigcap_n\brk{V_n}$, where $V_n=\{\sigma\in\strings: (n,\sigma)\in V\}$ is  $\Delta^1_1$. 
\end{proposition}

A set $X\in\cantor$ is $\Delta^1_1$-dominated if for any function $f\leq_hX$, there is a $\Delta^1_1$ function $g$ such that $\exists m\forall n>m(g(n)>f(n))$, we say $g$ dominates $f$. 
Recall that  $\omega_1^{X}$ denotes the first ordinal that is not recursive in $X$, and $\omega_1^{ck}$ denotes $\omega_1^{\emptyset}$.
There is a deep connection between $\Pi^1_1$ random and $\Delta_1^1$ random as displayed below.

\begin{theorem}[Kjos-Hanssen,Nies,Stephan and Yu \cite{kjos2010higher}]
	\label{thm:about-pi11-random}
	For any $X\in\cantor$, the followings are equivalent:
	\begin{enumerate}
		\item $X$ is $\Pi^1_1$-random.
		\item $X$ is $\Delta^1_1$-random and $\omega_1^X=\omega_1^{ck}$.
		\item $X$ is $\Delta^1_1$-random and $\Delta^1_1$-dominated, i.e. for any function $f\leq_h X$, there is a $\Delta^1_1$ function that dominates $f$.
	\end{enumerate}
\end{theorem}

The fact that $\Pi^1_1$ random is always $\Delta^1_1$ dominated actually enables us to convert a hyper reduction to a uniformly total reduction in some sense. That is, the following theorem of C.T.Chong and Yu holds. 

\begin{theorem}[C.T.Chong and L.Yu \cite{chong2015randomness}]
		\label{thm:pi11-random-comptation}
	For any $\Pi^1_1$ random real $X$ and $Y\leq_h X$, there is a notation $a\in\mathcal O$, a function $f\leq_TH_a$ and an oracle functional $\Phi$, such that for every $n$, $$Y(n)=\Phi^{X\oplus H_a\upto f(n)}(n)[f(n)]\downarrow,$$
	where we use $H_a$ to denote the H-sets upto $a\in\mathcal O$, and for any oracle $Z$, $H_b^Z$ is the H-sets relative to $Z$ upto some $b\in\mathcal O^Z$.
\end{theorem}

This theorem led to a higher version of Demuth's theorem in the paper of C.T.Chong and Yu \cite{chong2015randomness}. Likewise, in this paper, the above theorem will play an important role when we try to characterize the Hausdorff dimension of a maximal antichain of hyperdegrees.

There is a notion in classical recursion theory that's similar to be $\Delta^1_1$-dominated, i.e. the hyperimmune-free notion. We say a set $A$ is hyperimmune-free if for any function $f\leq_T A$, there is a recursive function $g$ such that $\exists m\forall n>m(g(n)>f(n))$. That is to say, every function computable in $A$ is computably dominated. A degree is of hyperimmune-free degree if it contains a hyperimmune-free set.
There is a well known basis theorem that says the hyperimmune-free sets in a way are very common. 

\begin{theorem}[The hyperimmune-free basis theorem,Jockusch and Soare \cite{jockusch1972pi01}]
	\label{thm:HIP basis}
	Every non-empty $\Pi^0_1$ class contains a member of hyperimmune-free degree.
\end{theorem}

Similar to theorem \ref{thm:pi11-random-comptation}, being hyperimmune-free can convert a Turing reduction to a truth-table reduction. 

\begin{theorem}[Jockusch \cite{jockusch1969relationships}, Martin]
	\label{thm:hif-tt-reduction}
	For any set $A$, the following are equivalent:
	\begin{enumerate}
		\item $A$ is of hyperimmune-free degree;
		\item For all functions $f\leq_T A$, $f\leq_{tt}A$;
		\item For all sets $B\leq_T A$, $B\leq_{tt} A$.
	\end{enumerate}
\end{theorem}

we will need this theorem in section \ref{sec:antichain-turing}, when we discuss the Hausdorff dimension of any maximal antichain of Turing degrees. 

In section \ref{sec:chains-turing}, we need another notion that extends the hyperimmune-freeness, namely the recursively traceable sets. Firstly, we fix a recursive enumeration $n\mapsto D_n$ of the finite subsets of $\omega$, that is, not only the map is recursive, but $D_n$ is uniformly recursive for every $n$. A set $A\in\cantor$ is called recursively traceable if there is a computable function $b$ (called a bound), such that for every function $f\leq_T A$, there is a computable function $g$ (called a trace), so that $|D_{g(n)}|\leq b(n)$ and $f(n)\in D_{g(n)}$ for all $n\in\omega$. A degree is recursively traceable if it contains a recursively traceable set. Obviously, being recursively traceable implies to be hyperimmune-free.

It is shown by Terwijn and Zambella that the bound in the definition of the recursively traceable sets can be very arbitrary. 

\begin{theorem}[Terwijn and Zambella \cite{terwijn1997algorithmic}]
	\label{thm:trace-bound-arbitrary}
	A degree $\mathbf a$ is recursively traceable if and only if for any unbounded, non-decreasing computable function $p:\omega\to\omega$ s.t. $p(0)>0$, and any function $f\leq_T \mathbf a$, there is a recursive function $g$ such that \begin{enumerate}
		\item $|D_{g(n)}|\leq p(n)$ for all $n$;
		\item $f(n)\in D_{g(n)}$ for all $n$.
	\end{enumerate}
\end{theorem}

A set $A\in\cantor$ is of minimal degree if it is non-computable, and every set $A$ computed is either computable, or computes $A$. Given a countable class of reals $\mathcal A$, a set $B$ is called minimal cover of $\mathcal A$ if for any $A\in\mathcal A$, $A<_T B$ and that for any $C\leq_T B$, there exists an $A\in\mathcal A$ such that either $C\leq_T A$ or $A\oplus C\geq_T B$. By results of Sacks \cite{sacks1963degrees}, every countable class of reals has a minimal cover. 

The foundational framework employed in this study relies on different variants of effective fractal dimensions. We provide a concise introduction to these concepts here; comprehensive treatments can be found in Downey and Hirschfeldt's seminal work \emph{Algorithmic Randomness and Complexity} \cite{downey2010algorithmic}. Recall that an order function $h:\omega\to\omega$ is one that's computable, non-decreasing, unbounded, and $h(0)>0$. Given a supermartingale $d$, the $h$-success set of $d$ is denoted $$S_h[d]=\{X\in\cantor: \limsup_n\frac{d(X\upto n)}{h(n)}=\infty\}$$It is easy to check that the function $n\mapsto 2^{(1-s)n}$ is an order function for all $s\in\mathbb Q\cap[0,1]$. There is a well known characterization of the Hausdorff dimension in terms of supermartingales.

\begin{theorem}[J.H.Lutz \cite{lutz2000gales,lutz2003dimensions}]
	\label{thm:martingale-characterization-dimension}
	Given any class $\mathcal R\subseteq \cantor$, we have:$$\dim_H(\mathcal R)=\inf\{s\in\mathbb Q\cap[0,1]:\mathcal R\subseteq S_{2^{(1-s)n}}[d]\text{ for some supermartingale }d\}.$$
\end{theorem}

This theorem motivates Lutz to consider effective Hausdorff dimension. He define the effective Hausdorff dimension of a class $\mathcal R\subseteq \cantor$ by $$\dim(\mathcal R)=\inf\{s\in\mathbb Q\cap [0,1]: \mathcal R\subseteq S_{2^{(1-s)n}}[d]\text{ for some r.e. supermartingale }d\}.$$ For a real $A\in\cantor$, the effective Hausdorff dimension of $A$ is $\dim(A)=\dim(\{A\})$.

Since there is an optimal r.e. supermartingale $d:\strings\to\mathbb R^{>0}$, i.e. one such that for any r.e. supermartingale $d'$, there is a constant $c\in\mathbb Q^{> 0}$ such that $d(\sigma)\geq cd'(\sigma)$ for all $\sigma\in\strings$, we conclude that $\dim(\mathcal A)=\inf\{s\in\mathbb Q\cap[0,1]:\mathcal A\subseteq S_{2^{(1-s)n}}[d]\}$. In particular, this implies the following.

\begin{proposition}[folklore]
	\label{prop:dim(R)=sup(dim(A))}
	Given any class $\mathcal R\subseteq\cantor$, we have $$\dim(\mathcal R)=\sup\{\dim(A):A\in \mathcal R\}.$$
\end{proposition}

\begin{proof}
	Given any class $\mathcal R$, by the discussion above, fix an optimal supermartingale $d$. Fix some $s\in\mathbb Q$ such that $\mathcal R\subseteq S_{2^{(1-s)n}}[d]$, then $A\in S_{2^{(1-s)n}}[d]$ for every $A\in\mathcal R$. So we have $\dim(A)\leq \dim(\mathcal R)$ for all $A\in\mathcal R$, that is $$\sup\{\dim(A):A\in\mathcal R\}\leq \dim(\mathcal R)$$
	Suppose the above inequation is strict. Fix an $r\in\mathbb Q$ such that $$\sup\{\dim(A):A\in\mathcal R\}<r\leq \dim(\mathcal R)$$ 
	By the definition, $\mathcal R\not\subseteq S_{2^{(1-r)n}}[d]$, thus there is some $B\in \mathcal R$ such that $B\notin S_{2^{(1-r)n}}[d]$, this implies $\dim(B)\geq r$, a contradiction.
\end{proof}

There is a beautiful characterization of the effective Hausdorff dimension in terms of Kolmogorov complexity, see Mayordomo  \cite{mayordomo2002kolmogorov}. That is, for any real $A\in\cantor$, we have $$\dim(A)=\liminf_{n\to\infty}\frac{K(A\upto n)}{n}=\liminf_{n\to\infty}\frac{C(A\upto n)}{n}.$$

Consequently, the effective Hausdorff dimension of a set is mathematically equivalent to the minimal asymptotic growth rate of its Kolmogorov complexity. This characterization, derived from the theoretical framework of Kolmogorov complexity, yields the following fundamental proposition that will serve as a cornerstone throughout our paper.

\begin{proposition}[Folklore]
	\label{thm:dim-deg}
	Let $A\leq_mB$, where $\leq_m$ is the many-one reduction. There is a $C\equiv_m B$ such that $A$ and $C$ have the same effective Hausdorff dimension. The same result holds for any weaker reduction, such as Turing reduction and hyperarithmetic reduction.
\end{proposition}

\begin{proof}
	Fix a fast-growing computable function $f:\omega\to\omega$ (such as the Ackerman function). And for each $n$, replace $A(f(n))$ to be $B(n)$, and call the resulting set $C$. Then apparently $C\equiv_m B$. And $C$ have the same effective Hausdorff dimension with $A$. Since the value of $K(A\upto n)/n$ and $K(C\upto n)/n$ are very similar for every $n$.
\end{proof}

We remark that this theorem can be partially relativized in the sense that given any oracle $Z$, if $A\leq_m B$, then there is some $C\equiv_m B$ such that $A$ and $C$ have the same effective-in-$Z$ Hausdorff dimension. The same holds for weaker reduction such as Turing reduction and Hyperarithmetic reduction.

We say a class $\mathcal A\subseteq \cantor$ is a class of Turing degrees if $\forall X\in \mathcal A\forall Y\in \cantor(X\equiv_T Y\to Y\in\mathcal A)$. The above proposition shows that we can identify a class of degrees with its downward closure in terms of Hausdorff dimension.

\begin{proposition}
	Given any class $\mathcal A$ of Turing degrees, let $ D_T(\mathcal A)=\{Y\in\cantor:\exists X\in \mathcal A(Y\leq_T X)\}$. Then we have $$\dim_{H}(\mathcal A)=\dim_{H}(D_T(\mathcal A)).$$ Let $U_T(\mathcal A)=\{Y\in\cantor:\exists X\in\mathcal A(X\leq_TY)\}$, then we have $$\dim_H(\mathcal A)\leq \dim_H(U_T(\mathcal A)).$$ And the above results also holds when "Turing reduction" is replaced by "hyperarithmetic reduction".
\end{proposition}

\begin{proof}
	The latter is obvious, since $\mathcal A\subseteq U_T(\mathcal A)$. To prove the former, notice that we only need to show that for any oracle $Z\in\cantor$, $$\dim^Z(\{X\in\cantor:X\equiv_T A\})=\dim^Z(\{X\in\cantor:X\leq_T A\})$$ for any $A\in\cantor$. But this is followed by the partial relativization of proposition \ref{thm:dim-deg}.
\end{proof}

Recall that we say a function $f:\omega\to\omega$ is DNR if $f(e)\not=\Phi_e(e)$ for any $e\in\omega$. A set is of DNR degree if it computes a DNR function. There is a deep connection between the notion of DNR and effective Hausdorff dimension. 

\begin{proposition}[Terwijn \cite{terwijn2004complexity}]
	\label{prop:DNR-and-effective-dimension}
	Given any set $A\in\cantor$, if $\dim(A)>0$, then $A$ is of DNR degree.
\end{proposition}

In section \ref{sec:antichain-turing}, we also care for the packing dimension of a maximal antichain of Turing degrees. In this case, we also have a theorem similar to theorem \ref{thm:martingale-characterization-dimension}. Recall that we say a supermartingale $d$ strongly $s$-succeed on a set $A$ if $\liminf_{n\to\infty}d(A\upto n)/2^{(1-s)n}=\infty$.
It is shown by Athreya, Hitchcock, Lutz, and Mayordomo \cite{athreya2007effective} that the packing dimension of a class $\mathcal A$ is equal to the infimum of the rational $s$ such that there is some supermartingale $d$ strongly $s$-succeed on all $A\in\mathcal A$.

So the effective packing dimension of a class $\mathcal R$ is defined by $$\Dim(\mathcal R)=\inf\{s\in\mathbb Q:\text{ some r.e. super-martingale s-succeeds strongly on all }A\in\mathcal R\}.$$And for a single set $A\in\cantor$, its effective packing dimension is $\Dim(A)=\Dim(\{A\})$.

The characterization of the effective packing dimension of a real in terms of the Kolmogorov complexity shows that $$\Dim(A)=\limsup_{n\to\infty}\frac{K(A\upto n)}{n}=\limsup_{n\to\infty}\frac{C(A\upto n)}{n}.$$

There are also propositions for effective packing dimension that's similar to proposition \ref{thm:dim-deg} and  \ref{prop:dim(R)=sup(dim(A))}. Since every supermartingale is r.e. in some real, we have the following two equations:\begin{enumerate}
	\item $\dim_H(\mathcal A)=\inf\{\dim^Z(\mathcal A):Z\in\cantor\}$,
	\item $\dim_P(\mathcal A)=\inf\{\Dim^Z(\mathcal A):Z\in\cantor\}$.
\end{enumerate}
Where $\dim^Z,\Dim^Z$ denotes the effective-in-$Z$ Hausdorff and packing dimension separately. 

The remainder of this section is devoted to determining the Hausdorff dimension of the two antichains introduced by C.T. Chong and L. Yu in their paper \cite{chong2016measure}.

Fix $\mathcal A$ as the antichain of Turing degrees constructed in theorem \ref{thm:Yu's-maximal-antichain}. And fix an oracle $Z\in\cantor$. We will analyze the effective Hausdorff dimensions of $\mathcal A$ relative to $Z$.

Firstly, pick a $\Pi^1_1(Z)$-random set $R$, it is followed by definition that $R$ is also $\Pi^1_1$-random. By theorem \ref{thm:Yu's-maximal-antichain}, fix an $X\in \mathcal A$ such that $R<_T X$. Since $R$ is $\Pi^1_1(Z)$-random, it is also 1-$Z$-random. So it must be that $$\dim^Z(R)=1,$$ Thus by theorem \ref{thm:dim-deg} relative to $Z$, there is a set $Y$ such that $Y\equiv_T X$ and $\dim^Z(Y)=\dim^Z(R)=1$, thus $$\dim^Z(\deg_T(X))=1.$$where $\deg_T(X)$ denotes $\{Y\in\cantor:Y\equiv_T X\}$. Hence by corollary \ref{prop:dim(R)=sup(dim(A))}, the effective-in-$Z$ Hausdorff dimension of $\mathcal A$ is also 1.

Overall, we have proved that the Hausdorff dimension of class $\mathcal A$ is one. The same technique can be used to show that the Hausdorff dimension of the antichain of hyperdegrees constructed in theorem \ref{thm:Yu's-maximal-antichain} is also 1. Hence we have proved the following result.

\begin{theorem}
	The antichains of Turing and hyperdegrees in theorem \ref{thm:Yu's-maximal-antichain} both have Hausdorff dimension 1.
\end{theorem}

It is this very observation that motivates the subject of this paper. We are thus led to ask: Does every maximal antichain of Turing degrees have Hausdorff dimension 1? If not, how small can its dimension be? More concretely, is it possible to construct a maximal antichain of Turing degrees with Hausdorff dimension strictly less than 1? Furthermore, what happens when we consider chains of Turing degrees instead of antichains, or when we examine antichains of hyperdegrees? This paper provides definitive answers to several of these questions.

\section{Antichains In Hyperdegrees}
\label{sec:antichain-hyper}

In this section, we study the Hausdorff dimension of antichains of hyperdegrees. By applying the $\Delta^1_1$-domination property of $\Pi^1_1$-random reals---as stated in Theorem \ref{thm:pi11-random-comptation}---we firstly show that, relative to any oracle $Z$, any non-HYP real $X$ hyperarithmetically-reducible to a $\Pi^1_1(Z)$-random real $R$ is essentially of $\Pi^1_1(Z)$-random hyperdegree.

\begin{lemma}
	\label{thm:pi11(Z)-computation}
	Given any oracle $Z\in\cantor$. If real $R$ is $\Pi^1_1(Z)$-random, and $\emptyset<_hX\leq_h R$, then there exists a real $Y\in\cantor$ such that $Y\equiv_h X$ and $Y$ is also $\Pi^1_1(Z)$-random.
\end{lemma}

Note that this is not just a simple relativization of the higher Demuth's theorem presented in \cite{chong2015randomness}, but rather a partial one.

\begin{proof}
	Fix an oracle $Z\in\cantor$, a $\Pi^1_1(Z)$-random real $R$, and a non-HYP real $X\leq_h R$. 
	
	Since $R$ is $\Pi^1_1(Z)$-random, it is also $\Pi^1_1$-random, and $X\leq_h R$. So by lemma \ref{thm:pi11-random-comptation}, fix a notation $a\in\mathcal O$, a function $f\leq_T H_a$ and an oracle functional $\Phi$ that satisfy $$X(n)=\Phi^{R\oplus H_a\upto f(n)}(n)[f(n)]\downarrow.$$ Assume without loss of generality that $f$ is non-decreasing. And the functional $\Phi$ satisfies that for all $n\in\omega$, $\tau\in 2^{f(n)}$, we have $\Phi^{\tau\oplus H_a\upto f(n)}(n)[f(n)]\downarrow\in\{0,1\}$. Since otherwise, one can easily define another functional $\Psi$ that does the work and still satisfies $\Psi^{R\oplus H_a\upto f(n)}(n)[f(n)]\downarrow=X(n)$ for all $n$, simply let it be 0 if the computation doesn't halt in $f(n)$ steps or it halts but does not output 0 or 1.
	
	Now we follow the idea of Demuth, define a premeasure $\rho$, i.e. a function from $\strings$ to $\R^{\geq 0}$ that satisfies $\rho(\lambda)=1$ and $\rho(\sigma0)+\rho(\sigma1)=\rho(\sigma)$ for all $\sigma$. The definition of $\rho$ is as follows:
	$$\rho(\sigma)=\mu\left(\bigcup\left\{\brk{\tau}:|\tau|=f(|\sigma|)\wedge(\forall n<|\sigma|)\;\Phi^{\tau\oplus H_a\upto_{f(|\sigma|)}}(n)[f(|\sigma|)]\downarrow=\sigma(n)\right\}\right).$$
	Apparently, $\rho$ is a $\Delta^1_1$ function, i.e. a function such that the set $$\{(\sigma,p,q):p,q\in\mathbb Q\wedge p<\rho(\sigma)< q\}$$is $\Delta^1_1$. Indeed, for any $\sigma$, the set $$D_\sigma=\{\tau\in2^{ f(|\sigma|)}:(\forall n<|\sigma|)\;\Phi^{\tau\oplus H_a\upto_{f(|\sigma|)}}(n)[f(|\sigma|)]\downarrow=\sigma(n)\}$$ is uniformly recursive in $H_a$, hence the value $\mu( \brk{D_\sigma})=\rho(\sigma)$ is uniformly recursive in $H_a$, which means that the set $$\{(\sigma,p,q):p,q\in\mathbb Q\wedge p<\mu(\brk{D_\sigma})<q\}$$is recursive in $H_a$. Further more, $\rho$ clearly satisfies $$\rho(\sigma)=\rho(\sigma 0)+\rho(\sigma1),$$ since by the assumption we made for the functional, for all $\sigma\in\strings$, we have $$ \brk{D_\sigma}= \brk{D_{\sigma0}}\sqcup \brk{D_{\sigma1}}.$$ where $\sqcup$ is the non-intersected union.
	The function $\rho$ induces an outer measure $\mu_\rho^*$ by letting $\mu_\rho^*(\mathcal D)$ to be the infimum of $\sum_{\sigma\in U}\rho(\sigma)$ over all $U\subseteq \strings$ such that $\mathcal D\subseteq  \brk{U}$. Then using Caratheodory's technique, we can extend $\mu_{\rho}^*$ to be a measure, denote $\nu=\mu_\rho^*$. 
	One can show that $X$ will be an $\nu$-$\Delta^1_1(Z)$-random, but here we adopt a more straightforward method.
	
	For every $\sigma,\tau\in\strings$, let $\sigma<_L\tau$ to mean that $$\exists k<\min(|\sigma|,|\tau|)[\sigma(k)<\tau(k)\wedge (\forall l<k)\; \sigma(l)=\tau(l)]$$ Now for every $\sigma$ we let:$$\begin{array}{l}
		l(\sigma)=\displaystyle\sum_{|\tau|=|\sigma|\wedge \tau<_L\sigma}\nu(\brk{\tau}),\\
		r(\sigma)=l(\sigma)+\nu(\brk{\tau}).
	\end{array}$$
	So the functions $r,l$ are both recursive in $H_a$ (in the sense mentioned before).
	let $l_n^X=l(X\upto n),r_n^X=r(X\upto n)$ for every $n$. It is easy to show that for all $n$, we have: $$l_n^X\leq l_{n+1}^X\leq r_{n+1}^X\leq r_n^X.$$
	Since $X\notin HYP$ by definition, it must be that $\lim_n(r_n^X-l_n^X)=0$. In fact, we have $$r_n^X-l_n^X=\rho(X\upto n)=\mu(\brk{D_{X\upto n}}),$$ and if $\lim_n(r_n^X-l_n^X)>0$, it will be that $$\mu\left(\bigcap_n\brk{D_{X\upto n}}\right)>0.$$ and the class $\bigcap_n\brk{D_{X\upto n}}$ is essentially $\{Z\in\cantor: \forall n\;\Phi^{Z\oplus H_a\upto_{f(n)}}(n)[f(n)]\downarrow=X(n)\}$. By the Lebesgue density theorem, there is a set $A\in\cantor$ such that $$\lim_{m\to\infty}2^m\mu\left(\brk{A\upto m}\cap\bigcap_n \brk{D_{X\upto n}}\right)=1.$$Hence, let $m\in\omega$ be such that $$2^m\mu\left(\brk{A\upto m}\cap\bigcap_n\brk{D_{X\upto n}}\right)>\frac3 4.$$Let $\sigma=A\upto m$. Then the above means that $$2^m\mu\left(\left\{Z\in\cantor:Z\succ \sigma\wedge \forall n\; \Phi^{Z\oplus H_a\upto_{f(n)}}(n)[f(n)]\downarrow=X(n)\right\}\right)>\frac3 4.$$ In this way, $X$ will be $H_a$-computable in the following way: enumerate strings $\tau\succeq \sigma$ as $\{\tau_s:s\in\omega\}$ such that every $\tau\succeq \sigma$ will show up infinitely in this sequence; to compute $X(n)$, define two sets $L,R$, at step $s$, if $\Phi^{\tau_s\oplus H_a\upto_{f(n)}}(n)[f(n)]\downarrow=0$, put $\tau_s$ into $L$, if $\Phi^{\tau_s\oplus H_a\upto_{f(n)}}(n)[f(n)]\downarrow=1$, put it into $R$. And if at any step $s$ we find that $2^m\mu(\brk{L_s})>3/4$, halt and output 0 ($X(n)$ will be 0), if we find $2^m\mu( \brk{R_s})>3/4$, halt and output 1 ($X(n)$ will be 1). This is a contradiction since $X$ is not HYP.
	
	So there must be a real number $y\in\mathbb R$ such that $\{y\}=\bigcap_n[l_n^X,r_n^X]$. Let the $Y\in\cantor$ to be its binary expansion sequence. Then one can show that $Y\equiv_h X$. Indeed, given a string $\sigma$, let $0.\sigma$ denote the dyadic rational number that $\sigma$ represents, then $$\sigma\prec Y\Eqvl \exists k(l_k^X\leq 0.\sigma\leq r_k^X\wedge r_k^X-l_k^X< 2^{-|\sigma|}).$$We know that sequences $\{l_k^X\},\{r_k^X\}$ is recursive in $X\oplus H_a$, and $\lim_n(r_n^X-l_n^X)=0$, so $Y$ is recursive in $X\oplus H_a$. On the other hand, to compute $X$, assume if we have computed $\sigma_n=X\upto n$, using $y$ , we search for an $i<2$ such that $y\notin[l(\sigma_n i), r(\sigma_n i)]$, if found, $\sigma_{n+1}=\sigma_n(1-i)$, since the function $l,r$ is recursive in $H_a$, $X$ is computable from $y$ and $H_a$.
	
	But we have $Y\equiv_h X\leq_h R$, and $\omega_1^{R\oplus Z}=\omega_1^Z$, so $\omega_1^{Y\oplus Z}=\omega_1^Z$ as well. If we want to show that $Y$ is $\Pi^1_1(Z)$-random, we only have to show it is $\Delta^1_1(Z)$-random, by theorem \ref{thm:about-pi11-random} relative to $Z$. 
	
	Assume $Y$ is not $\Delta^1_1(Z)$-random for a contradiction, then by the property of the $\Delta^1_1(Z)$-randomness in fact  \ref{prop:delta11-random-vs-schnorr}, there is a notation $b\in\mathcal O^Z$, and an $H_b^Z$-Schnorr-test $V\subseteq\omega\times\strings$ such that $Y\in \bigcap_n \brk{V_n}$. Assume without loss of generality that $|b|_Z>|a|$. For a string $\sigma$, define $(\sigma)_\mu=[p,q]$ such that $p=\sum_{|\tau|=|\sigma|\wedge \tau<_L\sigma}2^{-|\tau|}$, while $q=p+2^{-|\sigma|}$. It's not hard to see that the real numbers in the interval $(\sigma)_\mu$ are exactly those with binary expansion sequences in $ \brk{\sigma}$, i.e. $$A\in \brk{\sigma}\Eqvl 0.A\in (\sigma)_\mu,$$ if we let $0.A$ to denote the real number with binary expansion sequence $A\in\cantor$. 
	
	Now define: $$\begin{array}{c}
		\hat{V}_n=\{\sigma\in\strings:\exists \tau\exists k(\tau \text{ is the }k\text{-th string enumerated in }V_n\;\wedge\\
		\;\exists p,q\in\mathbb Q([p,q]=(\tau)_\mu\wedge [l(\sigma),r(\sigma)]\subseteq [p-2^{-(n+k+3)},q+2^{-(n+k+3)}]))\}.
	\end{array}$$
	Note that the $\{\hat{V}_n\}$ is uniformly recursively enumerable in $H_{b}^Z$. Since $Y\in \brk{ V_n}$, we can prove that $X\in \brk{\hat{V}_n}$. Fix the $k$-th string $\tau$ in $V_n$ such that $\tau\prec Y$, then $y\in (\tau)_\mu$, let $(\tau)_\mu=[p,q]$. Since $\nu(\brk{X\upto m})=r_m^X-l_m^X\to 0$ when $m\to \infty$, let $m$ be sufficiently large such that $\nu(\brk{X\upto m})<2^{-(n+k+3)}$. By the definition of $y$, we have $y\in[l(X\upto m),r(X\upto m)]$, so $$\begin{aligned}
		\relax[l(X\upto m),r(X\upto m)]&\subseteq [y-2^{-(n+k+3)},y+2^{-(n+k+3)}]\\
		&\subseteq [p-2^{-(n+k+3)},q+2^{-(n+k+3)}].
	\end{aligned}$$ Hence the string $X\upto m$ will eventually be enumerated into $\hat{V}_n$. Thus we have $X\in  \brk{\hat{V}_n}$ for all $n$. 
	
	Further define $$U_n=\{\tau\in\strings:\exists\sigma\in \hat{V}_n(|\tau|=f(|\sigma|)\wedge (\forall k<|\sigma|)\;\Phi^{\tau\oplus H_a\upto _{f(|\sigma|)}}(k )[f(|\sigma|)]=\sigma(k))\}.$$
	Then $\{U_n\}$ is uniformly recursively enumerable in $H_{b}^Z$, and obviously $R\in \bigcap_n \brk{U_n}$, since $X\in \brk{\hat{V}_n}$ for all $n$. But we can write $$U_n=\bigcup_{\sigma\in \hat{V}_n}D_{\sigma}.$$Note that for $\sigma\in \hat{V}_n$, $\mu(\brk{ D_{\sigma}})=\nu( \brk{\sigma})$. And for incomparable strings $\sigma,\tau\in \hat{V}_n$, $D_\sigma$ and $D_\tau$ clearly don't intersect. Next, by the definition of $\hat{V}_n$, for every $\sigma\in \hat{V}_n$, $$\nu( \brk{\sigma})=r(\sigma)-l(\sigma)\leq 2^{-(n+k+2)}+2^{-|\tau|}$$ for some $k$ and the $k$-th element $\tau\in V_n$.
	
	By the above analysis, we have the following estimate: $$\mu( \brk{U_n})\leq \mu( \brk{V_n})+\sum_{k=0}^\infty2^{-(n+k+2)}<2^{-n+1}.$$
	Hence $\{U_{n+1}\}$ is a $H_{b}^Z$-ML-test that covers $R$. This implies $R$ is not a $\Delta^1_1(Z)$-random, a contradiction.
\end{proof}

\begin{theorem}
	\label{thm:maximal-antichain-hyperdegrees}
	Every maximal antichain of hyperdegrees has Hausdorff dimension 1.
\end{theorem}

\begin{proof}
	Fix such a maximal antichain of hyperdegrees $\mathcal A$. Fix any oracle $Z$, we examine the effective in $Z$ Hausdorff dimension of $\mathcal A$.
	
	First of all, pick any $\Pi^1_1(Z)$-random real $R$. There are two cases to consider:
	
	Case 1: If $R\in \mathcal A$, then by the relativization of corollary \ref{prop:dim(R)=sup(dim(A))} to $Z$, we have $$\dim^Z(\mathcal A)\geq \dim^Z(R)=1.$$ 
	
	Case 2: If $R\notin \mathcal A$, now by the maximal property of $\mathcal A$, there is some $X\in \mathcal A$ such that $R\leq_h X$ or $X\leq_h R$. If $R\leq_hX$, by the relativized version of theorem \ref{thm:dim-deg}, the hyperdegree of $X$ has effective-in-$Z$ Hausdorff dimension 1. Hence $$\dim^Z(\mathcal A)\geq\dim^Z(\deg_h(X))=1.$$
	Lastly, if it is the case that $X\leq_h R$, by lemma \ref{thm:pi11(Z)-computation} just proved, there is some $Y\in\cantor$ such that $Y\equiv_h X$ and $Y$ is also a $\Pi^1_1(Z)$-random. Thus, we still have $\dim^Z(\mathcal A)\geq \dim^Z(\deg_h(X))=1$. 
\end{proof}

\section{Chains in Turing degrees}
\label{sec:chains-turing}

In this section, we consider the Hausdorff dimension of maximal chains of Turing degrees. We show that there exists a maximal $\omega_1$-chain of Turing degrees with Hausdorff dimension 0. Such a chain is constructed within the class of recursively traceable reals. Using a $\Pi^1_1$-inductive principle proposed by C.T. Chong and L. Yu in \cite{chong2009pi11}, we conclude that this chain can further be $\Pi^1_1$ under the assumption $\omega_1 = (\omega_1)^L$, where $L$ denotes G\"odel's constructible universe.

In recursion theory, it is known that every recursively traceable degree has effective Hausdorff dimension 0. This follows from the result of Kjos-Hanssen, Merkle, and Stephan \cite{kjos2011kolmogorov}, who showed that every recursively traceable set is neither DNR nor high. However, in this paper, we provide a direct proof of this fact.

\begin{lemma}
	\label{lem:traceable-dim-zero}
	Every recursively traceable set has effective Hausdorff dimension 0.
\end{lemma}

\begin{proof}
	Fix a set $A$ that is recursively traceable. This means, there exists a recursive function $h$ (called the bound), such that for any function $f\leq_TA$, there is a recursive function $g$ such that $|D_{g(n)}|\leq h(n)$ and $f(n)\in D_{g(n)}$ for all $n$.
	
	Now assume $\dim(A)>0$, fix a rational $s$ such that $\dim(A)>s>0$ and that $s^{-1}\in\mathbb N$. By definition, $A\notin S_{2^{(1-s)n}}[d]$, where $d$ is an optimal r.e. super-martingale. Let $p(n)=s^{-1}\cdot(n+\log h(n))$, so $p$ is recursive. And the map $n\mapsto A\upto p(n)$ is recursive in $A$. So by definition, there is a recursive function $g$ such that for every $n$, $|D_{g(n)}|\leq h(n)$ and $A\upto p(n)\in D_{g(n)}$, here we think of strings as their G\"odel codes.
	
	We can assume that inside every $D_{g(n)}$, there are only strings with length $p(n)$, this is legitimate since if not we can effective find a new $g^*$ satisfies this property. Given any string $\sigma$, we use $[\sigma]^{\preceq }$ to denote the set of finite strings that extends $\sigma$, i.e. $[\sigma]^\preceq=\{\tau\in\strings:\tau\succeq \sigma\}$. 
	
	Now we define a martingale $d_n$ for every $n$, such that:$$d_n(\sigma)=n\cdot2^{|\sigma|}\cdot\sum_{\tau\in D_{g(n)}\cap[\sigma]^{\preceq}}2^{-s|\tau|}.$$
	Hence for any $i<2$, we have:$$d_n(\sigma i)=n\cdot2^{|\sigma|+1}\cdot\sum_{\tau\in D_{g(n)}\cap[\sigma i]^{\preceq}}2^{-s|\tau|}.$$
	Because of the reason that $D_{g(n)}\cap[\sigma]^{\preceq}=(D_{g(n)}\cap[\sigma0]^{\preceq})\sqcup (D_{g(n)}\cap[\sigma1]^{\preceq})$, where $\sqcup$ means the non-intersected union, one can easily verify that every $d_n$ is an r.e. martingale. Let $\hat{d}=\sum_nd_n$, then $$\hat{d}(\lambda)=\sum_nd_n(\lambda)=\sum_n n\cdot\sum_{\tau\in D_{g(n)}}2^{-s|\tau|}\leq \sum_nn2^{-n}<\infty.$$So $\hat{d}$ is also an r.e. martingale. And if $\sigma\in D_{g(n)}$, then $$d_n(\sigma)=n\cdot2^{|\sigma|}\cdot\sum_{\tau\in D_{g(n)}\cap[\sigma]^{\preceq}}2^{-s|\sigma|}\geq n\cdot 2^{(1-s)|\sigma|}.$$ This implies $$\frac{d_n(\sigma)}{2^{(1-s)|\sigma|}}\geq n.$$ Since $A\upto p(n)\in D_{g(n)}$, so we actually have $$\limsup_n\frac{\hat{d}(A\upto n)}{2^{(1-s)n}}=\infty.$$
	that is $A\in S_{2^{(1-s)n}}[\hat{d}]$. But this implies $A\in S_{2^{(1-s)n}}[d]$, since $d$ is optimal. But this leads to a contradiction.
\end{proof}

As a corollary, the class $\{X\in\cantor:X\text{ is recursively traceable}\}$ has Hausdorff dimension zero. This is a straightforward deduction of lemma \ref{lem:traceable-dim-zero} and proposition \ref{prop:dim(R)=sup(dim(A))}.

\begin{corollary}
	\label{cor:traceable-dim-zero}
	The class $\{X\in\cantor:X\text{ is recursively traceable}\}$ has Hausdorff dimension zero.
\end{corollary}

Recall that given a countable class of reals $\mathcal A$, a minimal cover of $\mathcal A$ is a real $B$ such that for any $A\in\mathcal A$, $B>_T A$, and for any $C\leq_T B$ there exists some $ A\in \mathcal A$ such that either $C\leq_TA$ or $A\oplus C\geq_T B$. Sacks \cite{sacks1963degrees} showed that every countable class has a minimal cover using the notion of the now-called recursively pointed sacks forcing. The central lemma in this section extends Sacks' result in a sense that given any countable class of recursively traceable sets $\mathcal A$ that any two elements are Turing comparable, there is a recursively traceable set $B$ that forms a minimal cover for $\mathcal A$.

First of all, we introduce the notion of recursively pointed Sacks forcing. 

\begin{definition}
	\begin{enumerate}
		\item A function $T:\strings\to\strings$ is a (function) tree if $T(\sigma i)\succ T(\sigma)$ and $T(\sigma 0)\perp T(\sigma 1)$ for all $\sigma$ and $i<2$.
		\item $[T]$ denote the infinite paths of this tree, i.e. $$[T]=\{X\in\cantor: \exists Y\in\cantor\forall n(T(Y\upto n)\prec X)\}.$$
		\item For any two trees, we say $T\subseteq S$ if $[T]\subseteq [S]$.
		\item A tree $T$ is recursive if the map $T:\strings\to\strings$ is recursive. And it's recursively pointed if for all $X\in [T]$, $X\geq_T T$. It is uniformly pointed if there is an $e$ such that for every $X\in [T]$, $T=\Phi_e^X$.
		\item Given a function tree $T$ and a string $\sigma$, define another tree $\ext_\sigma(T)(\tau)=T(\sigma\tau)$ for all $\tau\in\strings$. 
	\end{enumerate}
\end{definition}

There are some fundamental facts about recursively pointed trees, which show that the notion of recursively pointed trees are in a way very robust.

\begin{lemma}[Sacks]
	\label{lem:Sacks-pointed-trees}
	Let $T$ be a recursively pointed function tree.
	\begin{enumerate}
		\item For every $Y\in\cantor$, $T\leq_TY$ if and only if $\exists X\in[T](X\equiv_T Y)$.
		\item For every $Y\in\cantor$, if $T\leq_T Y$, then there exists a recursively pointed tree $S$ such that $S\subseteq T$ and $S\equiv_T Y$.
		\item If $S$ is a function tree such that $S\subseteq T$ and $S\leq_T T$, then $S$ is also recursively pointed and $S\equiv_T T$.
	\end{enumerate}

\end{lemma}

Note that these results also hold when replace the "recursively pointed" to the "uniformly pointed". And the transition is uniform in the sense that, there is recursive functions $f,g$ such that if $T$ is uniformly pointed through $e_0$, and $T\leq_TY$ through $e_1$, then $\Phi_{f(e_0,e_1)}^Y$ computes a function tree $S$ such that $S\equiv_TY$ through $g(e_0,e_1)$.

\begin{proof}
	For 1, let $X$ be the unique path such that $\forall n(T(Y\upto n)\prec X)$. Then it is easy to show $X\oplus T\geq Y$,$X\leq_T T\oplus Y$, and since $X\geq_T T$, we have $X\geq_T Y$. Note that $X\leq_T Y$ is implied by $Y\geq _TT$.
	
	For 2, for every $\sigma$, define first that $\sigma\oplus Y$ to be a string with length $2|\sigma|$ and $$(\sigma\oplus Y)(2k)=\sigma(k),(\sigma\oplus Y)(2k+1)=Y(k)$$for all $k<|\sigma|$. Now define a function tree $S$ such that $S(\sigma)=T(\sigma\oplus Y)$ for every $\sigma$. It is easy to check that $S$ is a function tree and that $S\subseteq T$ and $S\oplus T\geq_T Y$, $S\leq_T T\oplus Y$. Combine the fact $Y\geq_T T$, and that the left most path of $S$ is recursive in $S$ and compute $T$, we have $S\equiv_T Y$. $S$ must be recursively pointed, since for any $X\in [S]$, $X\oplus T\geq Y$ hence $X\geq_T S$.
	
	For 3, $S$ is surely recursively pointed since the left most path of $S$ is recursive in $S$ and compute $T$.
	
	The uniformly pointed case is identical, one can analyze the above construction to see it.
\end{proof}

The following is central to this section. We use delicate uniformly pointed Sacks forcing and coding technique to extend the minimal covering theorem of Sacks. The technique used below essentially goes to C.T.Chong and L.Yu \cite{chong2007maximal}.

\begin{lemma}
	\label{lem:minimal-upper-bound}
	Given any countable sequence $\{A_i\}_{i\in\omega}$ that contains only recursively traceable sets, such that for any $i,j\in\omega$, either $A_i\leq_T A_j$ or $A_j\leq_T A_i$. Then for any real  $Z\in\cantor$, there exists a recursively traceable set $B$ that forms a minimal cover for sequence $\{A_i\}$ and satisfy $B''\geq_T Z$.
\end{lemma}

\begin{proof}
	By lemma \ref{lem:Sacks-pointed-trees}, given any uniformly pointed tree $T$ and a string $\sigma$, if $Y\geq_T T$, one can uniformly get a uniformly pointed tree $S\subseteq \ext_\sigma(T)$, such that $S\equiv_T Y$. Indeed, if $\Phi_e^X=T,\forall X\in[T]$, then there is a recursive $f:\omega\to\omega$ such that $\Phi_{f(\sigma)}^X=Y,\forall X\in[S]$, and a recursive $g$ such that $\Phi_{g(\sigma)}^X=S,\forall X\in[S]$. This procedure is called coding set $Y$ into tree $\ext_\sigma(T)$.
	
	By adding $\emptyset$ to the sequence, we can assume that $A_0\equiv_T \emptyset$. For any $i$, define $\hat{A}_i=\bigoplus_{k\leq i}A_k$, and since $\{A_i\}$ is a chain, $\hat{A}_i$ must be Turing equivalent to $A_k$ for some $k\leq i$. Thus every $\hat{A}_i$ is also recursively traceable, and $\hat{A}_i\leq_T \hat{A}_{i+1}$ for all $i\in\omega$.
	
	We will construct a sequence of uniformly pointed trees $\{T_n\}_{n\in\omega}$, a sequence $\langle e_n\rangle_{n\in\omega}$ such that $(\forall X\in [T_n])\;\Phi_{e_n}^X=T_n$, $\forall n(T_n\equiv_T\hat{A}_n\wedge T_n\supseteq T_{n+1})$, and the unique $B$ so that $\{B\}=\bigcap_n[T_n]$ will satisfy the whole requirements.
	
	At stage 0, we let $T_0=Id$, the identity tree, obviously $T_0\equiv_T \hat{A}_0\equiv_T\emptyset$. And fix $e_0$ to be a code for a Turing machine so that $\Phi_{e_0}^X=X$ for any $X\in\cantor$.
	
	At stage $s+1$, we split into the following steps to conduct:
	\begin{enumerate}
		\item[Step 1:](\emph{Coding set $Z$}) First of all, for every $k$, let $$\sigma_k=\left\{\begin{array}{ll}
			0^k1,&\text{ if }Z(s)=0,\\
			1^k0,&\text{ if }Z(s)=1.
		\end{array}\right.$$
		Let $S_k=\ext_{\sigma_k}(T_s)$. Then there is a recursive function $p$ such that $$\forall X\in[S_k]\; (\Phi_{p(k)}^X=S_k).$$ By the recursion theorem, fix a $k_0$ such that $$\forall X\in[S_{k_0}]\;(\Phi_{p(k_0)}^X=\Phi_{k_0}^X).$$Let $T_s^1=S_{k_0}$. Obviously, $T_s^1\equiv_T\hat{A}_s$ (although it maybe not uniform), since by inductive hypothesis, $T_s\equiv_T \hat{A}_s$.
		
		\item[Step 2:] (\emph{Coding set $\hat{A}_{s+1}$}) For every $k$, we code the set $\hat{A}_{s+1}$ into the tree $\ext_{0^k1}(T_s^1)$ to get $P_k$ by means mentioned in the beginning of the proof. So there are recursive functions $f,g$ such that $$(\forall X\in[P_k])\;[\Phi^X_{f(k)}=\hat{A}_{s+1}\wedge \Phi^X_{g(k)}=P_k].$$
		Using the recursion theorem, we can effectively find a $k_1$ such that $$(\forall X\in[P_{k_1}])\;[\Phi_{k_1}^X=\Phi_{g(k_1)}^X=S_{k_1}].$$
		Let $T_s^2=P_{k_1}$, then $T_s^2\equiv_T \hat{A}_{s+1}$, since $T_s^1\equiv_T\hat{A}_s$ and $\hat{A}_s\leq_T\hat{A}_{s+1}$, so $\hat{A}_{s+1}\geq_T T_s^2$. The other direction is obvious.
		
		\item[Step 3:] (\emph{forcing $B\not\leq_T \hat{A}_{i}$ for all $i\leq s$})  Choose a lexicographical least string $\sigma\in\strings$ such that for all $i,e\leq s$, $T_s^2(\sigma)\not\prec W_e^{T_i}$. Note $\sigma$ can be chosen so that $|\sigma|=(s+1)^2$. Then for any $k\in\omega$, define $R_s=\ext_{\sigma 0^k1}(T_s^2)$. So there is a recursive function $b$ such that $$(\forall X\in [R_k])\;\Phi^X_{b(k)}=R_k.$$ By recursion theorem, we effectively find a $k_2$ such that $$(\forall X\in[R_{k_2}])\;\Phi_{k_2}^X=\Phi_{b(k_2)}^X=R_{k_2}.$$ Let $T_s^3=R_{k_2}$, then $T_s^3\equiv_T \hat{A}_{s+1}$, since $T_s^3\subseteq T_s^2$ and $T_s^3\leq_T T_s^2$, and by lemma \ref{lem:Sacks-pointed-trees} (3).
		
		\item[Step 4:] (\emph{Forcing the minimal cover}) Let $e=s$. We split into the following cases to consider:
		\begin{enumerate}
			\item[Case 1:] If there exists a $\sigma$ and an $n$ such that $\forall \tau\succeq \sigma(\Phi_e^{T^3_s(\tau)})(n)\uparrow$, fix the least such $\sigma$ and $n$. For any $k$, we define $Q_k=\ext_{\sigma0^k1}(T_s^3)$. Then there exists a recursive function $q$ such that $$(\forall X\in[Q_k])\;\Phi_{q(k)}^X=Q_k.$$ By the recursion theorem again, we can effectively find a $k_3$ such that $$(\forall X\in[Q_{k_3}])[\Phi_{k_3}^{X}=\Phi_{q(k_3)}^X=Q_{k_3}].$$Similarly, Let $T_{s+1}=Q_{k_3},e_{s+1}=k_3$. We can also verify that $T_{s+1}\equiv_T \hat{A}_{s+1}$, since $T_{s+1}\equiv_T T_s^3\equiv_T \hat{A}_{s+1}$. In this case, if $B\in [T_{s+1}]$, then $\Phi_e^B(n)\uparrow$, so $\Phi_e^B$ is not a total function.
			
			\item[Case 2:] If case 1 fails, that is, for any $\sigma$ and $n$, $\exists \tau\succeq \sigma(\Phi_e^{T_s^3(\tau)}(n)\downarrow)$. At this point, we further split into two subcases to consider:
			\begin{enumerate}
				\item[Subcase 1:] If there is some $\sigma_0$ such that $$\forall n\forall \tau_0,\tau_1\succeq\sigma_0[\Phi_e^{T_s^3(\tau_0)}(n)\downarrow\wedge \Phi_e^{T_s^3(\tau_1)}(n)\downarrow\to \Phi_e^{T_s^3(\tau_0)}(n)=\Phi_e^{T_s^3(\tau_1)}(n)].$$ Again, for every $k$, define $Q_k=\ext_{\sigma0^k1}(T_s^3)$. Using the recursion theorem, we can find a $k_2$ such that $$(\forall X\in [Q_{k_3}])[\Phi_{k_3}^X=Q_{k_3}].$$
				Let $T_{s+1}=Q_{k_3}$ and $e_{s+1}=k_3$ then we are done. One can also check that $T_{s+1}\equiv_T\hat{A}_{s+1}$.
				
				\item[Subcase 2:] If Subcase 1 fails, then for all $\sigma$, there exist the least $n$ and least pair $(\tau_0,\tau_1)$ of incomparable strings extending $\sigma$ such that $\Phi_e^{T_s^3(\tau_0)}(n)\downarrow\not=\Phi_e^{T_s^3(\tau_1)}(n)\downarrow$. Similarly, define $Q_k=\ext_{0^k1}(T_s^3)$ for every $k$. And for every $k$, we can recursively-in-$Q_k$ find a subtree $H_k\subseteq Q_k$ such that for every $\tau$, $\Phi_e^{H_k(\tau)}(|\tau|)\downarrow$ and $(\exists n)\Phi_e^{H_k(\tau0)}(n)\downarrow\not=\Phi_e^{H_k(\tau1)}(n)\downarrow$. Just like how it is done in the classical construction of a $e$-splitting tree. Then there is also a recursive function $l$ such that $$(\forall X\in[H_{k}])[\Phi_{l(k)}^X=H_k]$$for all $k$. By recursion theorem again, find a $k_3$ such that $$(\forall X\in[H_{k_3}])[\Phi_{k_3}^X=\Phi_{l(k_3)}^X=H_{k_3}].$$ Let $T_{s+1}=H_{k_3}$ and $e_{s+1}=k_3$, it is routine to verify that $T_{s+1}\equiv_T\hat{A}_{s+1}$.
			\end{enumerate}
		\end{enumerate}
	\end{enumerate}
	
	Finally, let $B$ be the unique set that $\{B\}=\bigcap_n[T_n]$. Hence $B\geq_T\hat{A}_n$ for all $n$, since every $T_n$ is uniformly pointed and $T_n\equiv_T \hat{A}_{n}$.
	If $B\leq_T \hat{A}_{s}$ for some $s$, since $\hat{A}_s\equiv_T T_s$, assume $W_e^{T_s}=B$, then at stage $t+1=\max\{e,s\}$, step 3, there is some $\sigma$ such that $T_t^2(\sigma)\not\prec W_e^{T_s}$, and $T_{t+1}\subseteq \ext_{\sigma}(T_t^2)$, so $B\in[T_{t+1}]\subseteq [\ext_{\sigma}(T_s^2)]$, this leads to a contradiction.
	
	Assume $\Phi_e^B$ is total, then we are at stage $e+1$-case 2, if it subcase 1 holds, then $\Phi_e^B\leq_T \hat{A}_{e+1}$, otherwise it would be subcase 2 that holds and thus $\Phi_e^B\oplus \hat{A}_{e+1}\geq_T B$. 
	This proves that $B$ forms a minimal cover for the sequence $\{\hat{A}_n\}_{n\in\omega}$, and hence for $\{A_n\}$.
	
	Now we verify that $B''\geq_T Z$. Indeed, we are going to show that the sequence $\langle e_n\rangle_{n\in\omega}$ is recursive in $B''$. Obviously, $e_0$ can be recursively found. Assume inductively that $e_s$ was found, so that $\Phi_{e_s}^B=T_s$. Now we search for a $k_0$ such that either $T_s(0^{k_0}1)\prec B$ or $T_s(1^{k_0}0)\prec B$. If  the former holds, we know that $Z(s)=0$, and if it is the latter, $Z(s)=1$. No matter which one happens, $\Phi_{k_0}^B=T_s^1$. So next we search for a $k_1$ such that $T_s^1(0^{k_1}1)\prec B$, by definition, $\Phi_{k_1}^B=T_s^2$ and $\Phi_{f(k_1)}^B=A_{s+1}$. Further, we use the oracle $B''$ to find the least string $\sigma$ such that for any $i,e\leq s$, $T_s^2(\sigma)\not\prec W_e^{T_i}$, and search for a $k_2$ such that $T_s^2(\sigma 0^{k_2}1)\prec B$, by definition, $\Phi^B_{k_2}=T_s^3$. 
	
	To continue, we need the help of the oracle $B''$ to decide which case holds in step 4. If it is case 1, then use $B''$ to find the least $\sigma$ and least $n$ such that $\forall \tau\succeq \sigma(\Phi_s^{T_s^3(\tau)}(n)\uparrow)$. Further seek for a $k_3$ so that $T_s^3(\sigma0^{k_3}1)\prec B$. Thus, $\Phi_{k_3}^B=T_{s+1}$ and $e_{s+1}=k_3$, the induction can continue. If it is case 2, again we need $B''$ to decide which subcase it is. If subcase 1 holds, use $B''$ to find the least $\sigma$ such that for all $n$, all $\tau_0,\tau_1\succeq \sigma$, if $\Phi_s^{T_s^3(\tau_0)}(n)\downarrow$ and $\Phi_s^{T_s^3(\tau_1)}(n)\downarrow$, then $\Phi_s^{T_s^3(\tau_0)}(n)=\Phi_s^{T_s^3(\tau_1)}(n)$. Then seek for a $k_3$ such that $T_s^3(\sigma0^{k_3}1)\prec B$, and then $e_{s+1}=k_3$, the induction can continue. If subcase 2 holds, then one can just search for $k_3$ such that $T_s^3(0^{k_3}1)\prec B$, and so $e_{s+1}=k_3$, the induction can also continue.
	
	Last but not least, we check that $B$ is recursively traceable. Since every $\hat{A}_n$ is recursively traceable, by theorem \ref{thm:trace-bound-arbitrary} of Terwijn and Zambella, the bound can be very arbitrary. If we fix $b(n)=2^n$, then $b$ can be uniform bound function for every $\hat{A}_n$. Assume $\Phi_e^B$ is total, then at stage $e+1$, step 3, it must be case 2. So we split into two subcases. If subcase 1 holds, by definition we have $\Phi_e^B\leq_T\hat{A}_{e+1}$. Since $\hat{A}_{e+1}$ is recursively traceable, $\Phi_e^B$ must be traced by some recursive trace functions. If it is subcase 2, then by definition, we define a trace function $g\leq_T \hat{A}_{e+1}$ such that for any $n$, $g(n)$ is the canonical code for the finite set $\{\Phi_e^{T_{e+1}(\tau)}(|\tau|):\tau\in 2^n\}$. Thus, by the fact that $\hat{A}_{e+1}$ is recursively traceable, there is a recursive trace $h$ such that $$|D_{h(n)}|\leq b(n)\wedge g(n)\in D_{h(n)}$$ for all $n$. Now define a function $d$ such that for all $n$, $d(n)$ is the canonical code for the finite set $$F_n=\{m\in\omega: \exists a\in D_{h(n)}\exists k\leq b(n)[m\text{ is the }k\text{-th element of }D_a]\}.$$Hence, $d$ is recursive and $|D_{d(n)}|\leq b(n)^2$, $f(n)\in D_{g(n)}\subseteq D_{d(n)}$ for every $n$. So $B$ is a recursively traceable set with bound $b^2$.
\end{proof}

In particular, the above lemma produces a recursively traceable minimal cover for a single recursively traceable degree. 
Using these results, we can now construct a maximal chain in Turing degrees that is of Hausdorff dimension zero.

\begin{theorem}
	There is a maximal chain in Turing degrees that has Hausdorff dimension 0.
\end{theorem}

\begin{proof}
	Actually, we are going to construct a maximal $\omega_1$-chain that lays entirely in the class of recursively traceable sets, hence by lemma \ref{lem:traceable-dim-zero}, it has Hausdorff dimension zero.
	
	We inductively construct a recursively traceable degree $\mathbf d_\alpha$ for every ordinal $\alpha<\omega_1$, such that they satisfy:
	\begin{itemize}
		\item $\mathbf{d}_0$ is a minimal degree that's recursively traceable,
		\item $\forall \alpha<\omega_1(\mathbf d_\alpha<_T\mathbf d_{\alpha+1})$,
		\item For all $\alpha<\omega_1$, $\mathbf d_{\alpha+1}$ is a minimal cover above $\mathbf d_\alpha$.
		\item If $\eta<\omega_1$ is a limit ordinal, then $\mathbf d_\eta$ is the minimal cover for the countable set $\{\mathbf d_\alpha:\alpha<\eta\}$.
	\end{itemize}
	This will suffice to show that $\mathcal D=\{\mathbf d_\alpha:\alpha<\omega_1\}$ is a maximal chain. Since if degree $\mathbf a$ is not in $\mathcal D$, then firstly $\mathbf a$ cannot bound every $\mathbf d_\alpha$, since a degree can only bound countably many other degrees. So there is some $\alpha<\omega_1$ such that $\mathbf d_\alpha\not<_T \mathbf a$. Fix the minimal such $\alpha$, we have $\forall \beta<\alpha(\mathbf d_\beta<_T \mathbf a)$, but $\mathbf d_\alpha$ form a minimal upper bound for $\{\mathbf d_\beta:\beta<\alpha\}$, so $\mathbf a\not<_T \mathbf d_\alpha$. This proves the maximality of the chain $\mathcal D$ in the Turing degrees.
\end{proof}

In the paper of C.T.Chong and Yu \cite{chong2009pi11}, they proposed a $\Pi^1_1$ inductive principle $\mathfrak{J}$ as follows.

\begin{principle}[$\Pi^1_1$-inductive principle $\mathfrak{J}$]
	If a relation $P(X,Y)\subseteq \cantor\times\cantor$ is $\Pi^1_1$ and cofinally prograssive, i.e. $$\forall X\in\cantor\forall Z\in\cantor\exists Y\in\cantor(Y\geq_h Z\wedge P(X,Y))$$ Then there is a $\Pi^1_1$ class $\mathcal A\subseteq \mathcal C:=\{X\in\cantor:X\in L_{\omega_1^X}\}$ such that \begin{enumerate}
		\item $\ot(<_L\upto \mathcal A)=\omega_1$, i.e. the definable well order $<_L$ on the G\"odel's $L$ restricted on $\mathcal A$ has order type $\omega_1$.
		\item $\forall Y\in\mathcal A\exists X\in\mathcal C[X\text{  codes the class }\{Z:Z\in\mathcal A\wedge Z<_LY\}\wedge P(X,Y)]$.
	\end{enumerate}
	Here, $X$ codes a countable class $\mathcal D$ means $\mathcal D=\{(X)_n:n\in\omega\}$, where $m\in (X)_n\Eqvl \langle n,m\rangle\in X$. 
\end{principle}

C.T.Chong and Yu proved that $\omega_1=(\omega_1)^L$ holds if and only if the above principle holds. Here $L$ denotes the G\"odel's constructible universe.

\begin{lemma}[C.T.Chong and Yu \cite{chong2009pi11}]
	\label{lem:pi11-inductive-principle-holds}
	$\omega_1=(\omega_1)^L$ holds if and only if the $\Pi^1_1$ induction principle holds.
\end{lemma}

By this lemma, we can derive an existence corollary of a $\Pi^1_1$ maximal chain of Turing degrees that is of Hausdorff dimension zero, under the assumption that $\omega_1=(\omega_1)^L$.

\begin{theorem}
	Assume that $\omega_1=(\omega_1)^L$ holds, there is a $\Pi^1_1$ maximal chain of Turing degrees that's of Hausdorff dimension zero.
\end{theorem}

\begin{proof}
	Consider a relation $P(X,Y)$ such that $P(X,Y)$ holds if and only if: $$\begin{array}{c}
		\forall n[(X)_n\text{ is recursively traceable }]\wedge \forall n,m[(X)_n\leq_T(X)_m\vee (X)_m\leq_T(X)_n]\;\wedge\\ Y\text{ is recursively traceable }\wedge
		Y\text{ forms a minimal cover for the class }\{(X)_n:n\in\omega\}.
	\end{array}$$
	Obviously, this is a $\Pi^1_1$ relation, and by lemma \ref{lem:minimal-upper-bound}, the relation $P(X,Y)$ is cofinally progressive. Hence by lemma \ref{lem:pi11-inductive-principle-holds}, the $\Pi^1_1$ induction principle for $P(X,Y)$ holds. Thus there is a class $\hat{\mathcal A}\subseteq\mathcal C$ as prescribed. 
	
	Then we can use the well order $<_L\upto \hat{\mathcal A}$ to verify inductively that $\forall X,Y\in\hat{\mathcal A}[X\not=Y\to X<_T Y\vee Y<_T X]$. In fact, by the definition of $\hat{\mathcal A}$, given any two set $X,Y\in\hat{\mathcal A}$, assume $X<_L Y$, since $P(Z,Y)$ holds for some $Z\in\mathcal C$, such that $\{(Z)_n:n\in\omega\}=\{A\in \hat{\mathcal A}:A<_LY\}$, we have $Y>_T X$. 
	
	That is to say, if we let $\mathcal A=\{X:\exists Y\in\hat{\mathcal A}[Y\equiv_T X]\}$, then the class $\mathcal A$ forms a chain of Turing degrees. Further more, we can show it is maximal. Fix a set $Z\notin \mathcal A$, hence $Z$ is not the degree of any $X\in\hat{\mathcal A}$. And it is not the case that $Z>_T X$ for all $X\in\mathcal A$, since a degree can only bound countably many degrees. Fix an enumeration $(X_\alpha)_{\alpha\in\omega_1}$ for the class $\hat{\mathcal A}$ by order of $<_L$. Then there is some $\alpha<\omega_1$ such that $Z\not>_T X_\alpha$. Fix the minimal such $\alpha$, hence for any $\beta<\alpha$, we have $X_\beta<_T Z$. By the property of $\hat{\mathcal A}$, $X_\alpha$ forms a minimal upper bound for $\{X_\beta:\beta<\alpha\}$, hence $Z\not<_T X_\alpha$. This proves the maximality of the chain $\mathcal A$ in Turing degrees.
	
	Since $$X\in \mathcal A\Eqvl \exists e[\Phi_e^X\text{ is total }\wedge \Phi_e^X\in\hat{\mathcal A}\wedge X\leq_T\Phi_e^X].$$ The class $\mathcal A$ is also $\Pi^1_1$.
	And since every real in $\hat{\mathcal A}$ is recursively traceable, the class $\mathcal A$ has Hausdorff dimension zero.
\end{proof}

\section{Antichains in Turing degrees}
\label{sec:antichain-turing}

In this section, we explore the Hausdorff dimension of maximal antichains within the Turing degrees. We demonstrate that every maximal antichain of Turing degrees has effective Hausdorff dimension 1. However, we highlight that this result does not relativize to arbitrary oracles, as was done in Theorem \ref{thm:maximal-antichain-hyperdegrees}, with the exception of oracles of K-trivial degrees. This distinction underscores the unique properties of K-trivial degrees in the context of relativization.

The following lemma is the well known Demuth's theorem.

\begin{lemma}[Demuth's Theorem]
	Given any sets $X,Y\in\cantor$, if $X$ is 1-random, $Y$ is not computable and $Y\leq_{tt} X$, then there exists 1-random set $Z$ such that $Z\equiv_T Y$.
\end{lemma}

By the hyperimmune-free basis theorem \ref{thm:HIP basis}, there is a 1-random real that's hyperimmune-free, since there is obviously a $\Pi^0_1$ class that contains only 1-random reals. Combining with  Demuth's theorem, we have the following.

\begin{lemma}
	\label{lem:hif-random-computes-random}
	Given any real $X,R\in\cantor$, if $R$ is 1-random and hyperimmune-free, and $0<_T X\leq_T R$, then there is some real $Y\in\cantor$ such that $Y\equiv_T X$ and $Y$ is 1-random.
\end{lemma}

\begin{proof}
	Fix such $X,R$, since $R$ is hyperimmune-free and $X\leq_T R$, by theorem \ref{thm:hif-tt-reduction}, $X\leq_{tt}R$. Next, since $X\not\leq_T \emptyset$, by Demuth's theorem, there is some $Y\in\cantor$ such that $Y\equiv_T X$ and $Y$ is also 1-random.
\end{proof}

Using the above Lemma, we can show it is impossible for one to construct a maximal antichain of Turing degrees of effective Hausdorff dimension zero.

\begin{theorem}
	\label{thm:antichain-in-Turing-degrees}
	Every maximal antichain of Turing degrees has effective Hausdorff dimension one.
\end{theorem}

\begin{proof}
	Fix a maximal antichain of Turing degrees $\mathcal A$. Pick a 1-random real $R$ that is hyperimmune-free, if $R\in \mathcal A$, then by corollary \ref{prop:dim(R)=sup(dim(A))}, we have $$\dim(\mathcal A)\geq \dim(R)=1.$$
	Now assume $R\notin \mathcal A$, so there is some $X\in \mathcal A$ such that either $X\leq_T R$ or $R\leq_TX$. If it is the latter, then by theorem \ref{thm:dim-deg} and corollary \ref{prop:dim(R)=sup(dim(A))}, we have $$\dim(\mathcal A)\geq \dim(\deg_T(X))=1.$$ 
	If it is the former, i.e. $X\leq_T R$, since $X\in\mathcal A$, $X$ is not computable, so by lemma  \ref{lem:hif-random-computes-random}, there is some $Y\equiv_TX$ that is also 1-random. This implies $$\dim(\mathcal A)\geq \dim(\deg_T(X))=1.$$
	Thus, every situation implies $\mathcal A$ has effective Hausdorff dimension one.
\end{proof}

The above result relies heavily on Lemma \ref{lem:hif-random-computes-random}. However, the fully relativized version of this lemma states that, for any oracle $Z$, if a $1$-$Z$-random real $R$ is hyperimmune-free relative to $Z$, then for every $X \leq_T R\oplus Z$ with $X \not\leq_T Z$, there exists a $1$-$Z$-random real $Y$ such that $Y \oplus Z \equiv_T X \oplus Z$. Whether this fact can be partially relativized in the sense of Lemma \ref{thm:pi11(Z)-computation} is crucial for determining the full Hausdorff dimension of a maximal antichain of Turing degrees.

Unfortunately, partial relativization of this fact to any oracle is impossible. This impossibility stems from the result that the class of hyperimmune-free sets has Hausdorff dimension 0. This follows directly from Miller's result \cite{khan2017forcing}, which proves that all hyperimmune-free sets are non-DNR relative to $\emptyset'$, and thus by Proposition \ref{prop:DNR-and-effective-dimension}, they have effectively-in-$\emptyset'$ Hausdorff dimension 0.

Moreover, since every 2-random real is hyperimmune, one cannot rely on resolving this question by using stronger notions of randomness than 2-randomness. In fact, it is well-known that any 2-random real computes a 1-generic real, and the latter has effective Hausdorff dimension 0.

Nevertheless, we can relativize Theorem \ref{thm:antichain-in-Turing-degrees} to K-trivial oracles. This represents the furthest extent achievable with current techniques in this direction.

\begin{theorem}
	\label{thm:relative-antichain-turing}
	Assume $Z$ is K-trivial. Then every maximal antichain of Turing degrees has effective in $Z$ Hausdorff dimension 1.
\end{theorem}

\begin{proof}
	Fix a hyperimmune-free 1-random set $R$, since K-trivial is equivalent to low-for-1-randomness, $R$ is also 1-$Z$-random. 
	
	Now, for all $X\leq_T R$, since $R$ is hyperimmune-free, we have $X\leq_{tt} R$ by theorem \ref{thm:hif-tt-reduction}. Thus we can prove there is a 1-$Z$-random real $Y$ such that $Y\equiv_T X$, much like how we did it in theorem \ref{thm:maximal-antichain-hyperdegrees}. 
	
	Hence fix any maximal antichain $\mathcal A\subseteq \cantor$, if there is an $X\in\mathcal A$ such that $X\leq_T R$, then $X$ is Turing equivalent to a 1-$Z$-random set, and so the effective in $Z$ Hausdorff dimension of $\mathcal A$ is 1. Other situations are trivial.
\end{proof}

By contrast, when considering packing dimension, we obtain significantly more favorable results. Notably, a theorem of Downey and Greenberg establishes that the packing dimension of the class of minimal Turing degrees is exactly 1. 


\begin{theorem}[Downey and Greenberg \cite{downey2008turing}]
	\label{thm:minimal-packing-one}
	For every oracle $Z$, there exists a minimal degree $X$ such that the effective in $Z$ packing dimension of $X$ is one.
\end{theorem}

Using this theorem, it is easy to characterize the possible packing dimension of a maximal antichain of Turing degrees.

\begin{theorem}
	\label{thm:antichain-turing-packing}
	Every maximal antichain of Turing degrees has packing dimension one.
\end{theorem}

\begin{proof}
	Fix a maximal antichain $\mathcal A$ of Turing degrees, and an oracle $Z$. By theorem \ref{thm:minimal-packing-one}, there is a minimal degree $X$ such that $\Dim^Z(X)=1$. Hence, if $X\in \mathcal A$, then $\Dim^Z(\mathcal A)=1$. And if $X\notin \mathcal A$, there is some $Y\in\mathcal A$ such that either $X<_T Y$ or $Y<_TX$ by the maximality. 
	
	If it is the case that $X<_T Y$, then by lemma \ref{thm:dim-deg}, $\Dim^Z(\deg_T(Y))=1$, hence also $\Dim^Z(\mathcal A)=1$.
	
	It cannot be the case that $Y<_T X$, otherwise, since $Y\not\leq_T \emptyset$, and $X$ is of minimal degree, we would have $Y\geq_T X$, a contradiction.
\end{proof}

This result reminds us that the question of Hausdorff dimension of antichains of Turing degrees could be related to the question of Hausdorff dimension of the class of minimal degrees.
Define the class $$\mathrm{Min}=\{X\in\cantor:X\text{ is of minimal degree}\}$$ Then we have the following.

\begin{proposition}
	Suppose $\dim_H(\mathrm{Min})=r\leq 1$, then every maximal antichain of Turing degrees has Hausdorff dimension greater than or equal to $r$.
\end{proposition}

\begin{proof}
	The proof is identical with theorem \ref{thm:antichain-turing-packing}.
\end{proof}

However, the exact Hausdorff dimension of the class $\mathrm{Min}$ is still an open question in recursion theory. We even don't know if it has effective Hausdorff dimension 1. Several recursion theorists have been working on this subject, and have made some progresses, see Greenberg and Miller \cite{greenberg2011diagonally}, Khan and Miller \cite{khan2017forcing}, and Liu Lu \cite{liu2019dnr} for more information.

While substantial evidence suggests that the minimal attainable Hausdorff dimension for maximal antichains of Turing degrees is 1, countervailing observations complicate this hypothesis. Crucially, the class of K-trivial sets---known to be Turing incomplete \cite{downey2002trivial} and $\omega$-recursively enumerable \cite{nies2005lowness}---motivates a critical inquiry: Might there exist an oracle $Z\in\cantor$ (potentially even $\emptyset'$) capable of reducing the $Z$-effective Hausdorff dimension of certain maximal antichains below 1? This tension compels us to formalize the following open problem:

\begin{question}
	Is it true that given any oracle $Z$, any maximal antichain of Turing degrees has effective-in-$Z$ Hausdorff dimension 1? If not, in what complexity does the oracle have to be so that there is some maximal antichain of Turing degrees that has effective-in-$Z$ Hausdorff dimension $<1$? Can $\emptyset'$ do the job?
\end{question}


\bibliography{sn-bibliography}


\begin{thebibliography}{26}
\ifx \bisbn   \undefined \def \bisbn  #1{ISBN #1}\fi
\ifx \binits  \undefined \def \binits#1{#1}\fi
\ifx \bauthor  \undefined \def \bauthor#1{#1}\fi
\ifx \batitle  \undefined \def \batitle#1{#1}\fi
\ifx \bjtitle  \undefined \def \bjtitle#1{#1}\fi
\ifx \bvolume  \undefined \def \bvolume#1{\textbf{#1}}\fi
\ifx \byear  \undefined \def \byear#1{#1}\fi
\ifx \bissue  \undefined \def \bissue#1{#1}\fi
\ifx \bfpage  \undefined \def \bfpage#1{#1}\fi
\ifx \blpage  \undefined \def \blpage #1{#1}\fi
\ifx \burl  \undefined \def \burl#1{\textsf{#1}}\fi
\ifx \doiurl  \undefined \def \doiurl#1{\url{https://doi.org/#1}}\fi
\ifx \betal  \undefined \def \betal{\textit{et al.}}\fi
\ifx \binstitute  \undefined \def \binstitute#1{#1}\fi
\ifx \binstitutionaled  \undefined \def \binstitutionaled#1{#1}\fi
\ifx \bctitle  \undefined \def \bctitle#1{#1}\fi
\ifx \beditor  \undefined \def \beditor#1{#1}\fi
\ifx \bpublisher  \undefined \def \bpublisher#1{#1}\fi
\ifx \bbtitle  \undefined \def \bbtitle#1{#1}\fi
\ifx \bedition  \undefined \def \bedition#1{#1}\fi
\ifx \bseriesno  \undefined \def \bseriesno#1{#1}\fi
\ifx \blocation  \undefined \def \blocation#1{#1}\fi
\ifx \bsertitle  \undefined \def \bsertitle#1{#1}\fi
\ifx \bsnm \undefined \def \bsnm#1{#1}\fi
\ifx \bsuffix \undefined \def \bsuffix#1{#1}\fi
\ifx \bparticle \undefined \def \bparticle#1{#1}\fi
\ifx \barticle \undefined \def \barticle#1{#1}\fi
\bibcommenthead
\ifx \bconfdate \undefined \def \bconfdate #1{#1}\fi
\ifx \botherref \undefined \def \botherref #1{#1}\fi
\ifx \url \undefined \def \url#1{\textsf{#1}}\fi
\ifx \bchapter \undefined \def \bchapter#1{#1}\fi
\ifx \bbook \undefined \def \bbook#1{#1}\fi
\ifx \bcomment \undefined \def \bcomment#1{#1}\fi
\ifx \oauthor \undefined \def \oauthor#1{#1}\fi
\ifx \citeauthoryear \undefined \def \citeauthoryear#1{#1}\fi
\ifx \endbibitem  \undefined \def \endbibitem {}\fi
\ifx \bconflocation  \undefined \def \bconflocation#1{#1}\fi
\ifx \arxivurl  \undefined \def \arxivurl#1{\textsf{#1}}\fi
\csname PreBibitemsHook\endcsname

\bibitem[\protect\citeauthoryear{Yu}{2006}]{yu2006measure}
\begin{barticle}
\bauthor{\bsnm{Yu}, \binits{L.}}:
\batitle{Measure theory aspects of locally countable orderings}.
\bjtitle{The Journal of Symbolic Logic}
\bvolume{71}(\bissue{3}),
\bfpage{958}--\blpage{968}
(\byear{2006})
\end{barticle}
\endbibitem

\bibitem[\protect\citeauthoryear{Chong and Yu}{2016}]{chong2016measure}
\begin{barticle}
\bauthor{\bsnm{Chong}, \binits{C.}},
\bauthor{\bsnm{Yu}, \binits{L.}}:
\batitle{Measure-theoretic applications of higher demuth’s theorem}.
\bjtitle{Transactions of the American Mathematical Society}
\bvolume{368}(\bissue{11}),
\bfpage{8249}--\blpage{8265}
(\byear{2016})
\end{barticle}
\endbibitem

\bibitem[\protect\citeauthoryear{Chong and Yu}{2015}]{chong2015randomness}
\begin{barticle}
\bauthor{\bsnm{Chong}, \binits{C.T.}},
\bauthor{\bsnm{Yu}, \binits{L.}}:
\batitle{Randomness in the higher setting}.
\bjtitle{The Journal of Symbolic Logic}
\bvolume{80}(\bissue{4}),
\bfpage{1131}--\blpage{1148}
(\byear{2015})
\end{barticle}
\endbibitem

\bibitem[\protect\citeauthoryear{Downey and Greenberg}{2008}]{downey2008turing}
\begin{barticle}
\bauthor{\bsnm{Downey}, \binits{R.}},
\bauthor{\bsnm{Greenberg}, \binits{N.}}:
\batitle{Turing degrees of reals of positive effective packing dimension}.
\bjtitle{Information Processing Letters}
\bvolume{108}(\bissue{5}),
\bfpage{298}--\blpage{303}
(\byear{2008})
\end{barticle}
\endbibitem

\bibitem[\protect\citeauthoryear{Soare}{2016}]{soare2016turing}
\begin{botherref}
\oauthor{\bsnm{Soare}, \binits{R.I.}}:
Turing computability: theory and applications.
(No Title)
(2016)
\end{botherref}
\endbibitem

\bibitem[\protect\citeauthoryear{{Downey, Rodney G and Hirschfeldt, Denis
  R}}{2010}]{downey2010algorithmic}
\begin{bbook}
\bauthor{\bsnm{{Downey, Rodney G and Hirschfeldt, Denis R}}}:
\bbtitle{Algorithmic Randomness and Complexity}.
\bpublisher{Springer},
\blocation{New York}
(\byear{2010})
\end{bbook}
\endbibitem

\bibitem[\protect\citeauthoryear{Chong and Yu}{2015}]{chong2015recursion}
\begin{bbook}
\bauthor{\bsnm{Chong}, \binits{C.T.}},
\bauthor{\bsnm{Yu}, \binits{L.}}:
\bbtitle{Recursion Theory: Computational Aspects of Definability}
vol. \bseriesno{8}.
\bpublisher{Walter de Gruyter GmbH \& Co KG},
\blocation{Berlin}
(\byear{2015})
\end{bbook}
\endbibitem

\bibitem[\protect\citeauthoryear{Sacks}{1966}]{sacks1966degrees}
\begin{bbook}
\bauthor{\bsnm{Sacks}, \binits{G.E.}}:
\bbtitle{Degrees of Unsolvability}
vol. \bseriesno{55}.
\bpublisher{Princeton University Press},
\blocation{Princeton, New Jersey}
(\byear{1966})
\end{bbook}
\endbibitem

\bibitem[\protect\citeauthoryear{Kjos-Hanssen et~al.}{2010}]{kjos2010higher}
\begin{barticle}
\bauthor{\bsnm{Kjos-Hanssen}, \binits{B.}},
\bauthor{\bsnm{Nies}, \binits{A.}},
\bauthor{\bsnm{Stephan}, \binits{F.}},
\bauthor{\bsnm{Yu}, \binits{L.}}:
\batitle{Higher kurtz randomness}.
\bjtitle{Annals of Pure and Applied Logic}
\bvolume{161}(\bissue{10}),
\bfpage{1280}--\blpage{1290}
(\byear{2010})
\end{barticle}
\endbibitem

\bibitem[\protect\citeauthoryear{Jockusch and Soare}{1972}]{jockusch1972pi01}
\begin{barticle}
\bauthor{\bsnm{Jockusch}, \binits{C.G.}},
\bauthor{\bsnm{Soare}, \binits{R.I.}}:
\batitle{$\pi^0_1$ classes and degrees of theories}.
\bjtitle{Transactions of the American Mathematical Society}
\bvolume{173},
\bfpage{33}--\blpage{56}
(\byear{1972})
\end{barticle}
\endbibitem

\bibitem[\protect\citeauthoryear{Jockusch}{1969}]{jockusch1969relationships}
\begin{barticle}
\bauthor{\bsnm{Jockusch}, \binits{C.G.}}:
\batitle{Relationships between reducibilities}.
\bjtitle{Transactions of the American Mathematical Society}
\bvolume{142},
\bfpage{229}--\blpage{237}
(\byear{1969})
\end{barticle}
\endbibitem

\bibitem[\protect\citeauthoryear{Terwijn and
  Zambella}{1997}]{terwijn1997algorithmic}
\begin{botherref}
\oauthor{\bsnm{Terwijn}, \binits{S.A.}},
\oauthor{\bsnm{Zambella}, \binits{D.}}:
Algorithmic randomness and lowness
(1997)
\end{botherref}
\endbibitem

\bibitem[\protect\citeauthoryear{Sacks}{1963}]{sacks1963degrees}
\begin{barticle}
\bauthor{\bsnm{Sacks}, \binits{G.E.}}:
\batitle{On the degrees less than $0'$}.
\bjtitle{Annals of mathematics}
\bvolume{77}(\bissue{2}),
\bfpage{211}--\blpage{231}
(\byear{1963})
\end{barticle}
\endbibitem

\bibitem[\protect\citeauthoryear{Lutz}{2000}]{lutz2000gales}
\begin{bchapter}
\bauthor{\bsnm{Lutz}, \binits{J.H.}}:
\bctitle{Gales and the constructive dimension of individual sequences}.
In: \bbtitle{Automata, Languages and Programming: 27th International
  Colloquium, ICALP 2000 Geneva, Switzerland, July 9--15, 2000 Proceedings 27},
pp. \bfpage{902}--\blpage{913}
(\byear{2000}).
\bcomment{Springer}
\end{bchapter}
\endbibitem

\bibitem[\protect\citeauthoryear{Lutz}{2003}]{lutz2003dimensions}
\begin{barticle}
\bauthor{\bsnm{Lutz}, \binits{J.H.}}:
\batitle{The dimensions of individual strings and sequences}.
\bjtitle{Information and Computation}
\bvolume{187}(\bissue{1}),
\bfpage{49}--\blpage{79}
(\byear{2003})
\end{barticle}
\endbibitem

\bibitem[\protect\citeauthoryear{Mayordomo}{2002}]{mayordomo2002kolmogorov}
\begin{barticle}
\bauthor{\bsnm{Mayordomo}, \binits{E.}}:
\batitle{A kolmogorov complexity characterization of constructive hausdorff
  dimension}.
\bjtitle{Information Processing Letters}
\bvolume{84}(\bissue{1}),
\bfpage{1}--\blpage{3}
(\byear{2002})
\end{barticle}
\endbibitem

\bibitem[\protect\citeauthoryear{Terwijn}{2004}]{terwijn2004complexity}
\begin{bbook}
\bauthor{\bsnm{Terwijn}, \binits{S.A.}}:
\bbtitle{Complexity and Randomness}.
\bpublisher{na}, \blocation{???}
(\byear{2004})
\end{bbook}
\endbibitem

\bibitem[\protect\citeauthoryear{Athreya et~al.}{2007}]{athreya2007effective}
\begin{barticle}
\bauthor{\bsnm{Athreya}, \binits{K.B.}},
\bauthor{\bsnm{Hitchcock}, \binits{J.M.}},
\bauthor{\bsnm{Lutz}, \binits{J.H.}},
\bauthor{\bsnm{Mayordomo}, \binits{E.}}:
\batitle{Effective strong dimension in algorithmic information and
  computational complexity}.
\bjtitle{SIAM journal on computing}
\bvolume{37}(\bissue{3}),
\bfpage{671}--\blpage{705}
(\byear{2007})
\end{barticle}
\endbibitem

\bibitem[\protect\citeauthoryear{Chong and Yu}{2009}]{chong2009pi11}
\begin{barticle}
\bauthor{\bsnm{Chong}, \binits{C.}},
\bauthor{\bsnm{Yu}, \binits{L.}}:
\batitle{A $\pi$$^1_1$-uniformization principle for reals}.
\bjtitle{Transactions of the American Mathematical Society}
\bvolume{361}(\bissue{8}),
\bfpage{4233}--\blpage{4245}
(\byear{2009})
\end{barticle}
\endbibitem

\bibitem[\protect\citeauthoryear{Kjos-Hanssen
  et~al.}{2011}]{kjos2011kolmogorov}
\begin{barticle}
\bauthor{\bsnm{Kjos-Hanssen}, \binits{B.}},
\bauthor{\bsnm{Merkle}, \binits{W.}},
\bauthor{\bsnm{Stephan}, \binits{F.}}:
\batitle{Kolmogorov complexity and the recursion theorem}.
\bjtitle{Transactions of the American Mathematical Society}
\bvolume{363}(\bissue{10}),
\bfpage{5465}--\blpage{5480}
(\byear{2011})
\end{barticle}
\endbibitem

\bibitem[\protect\citeauthoryear{Chong and Yu}{2007}]{chong2007maximal}
\begin{barticle}
\bauthor{\bsnm{Chong}, \binits{C.T.}},
\bauthor{\bsnm{Yu}, \binits{L.}}:
\batitle{Maximal chains in the turing degrees}.
\bjtitle{The Journal of Symbolic Logic}
\bvolume{72}(\bissue{4}),
\bfpage{1219}--\blpage{1227}
(\byear{2007})
\end{barticle}
\endbibitem

\bibitem[\protect\citeauthoryear{Khan and Miller}{2017}]{khan2017forcing}
\begin{barticle}
\bauthor{\bsnm{Khan}, \binits{M.}},
\bauthor{\bsnm{Miller}, \binits{J.S.}}:
\batitle{Forcing with bushy trees}.
\bjtitle{Bulletin of Symbolic Logic}
\bvolume{23}(\bissue{2}),
\bfpage{160}--\blpage{180}
(\byear{2017})
\end{barticle}
\endbibitem

\bibitem[\protect\citeauthoryear{Greenberg and
  Miller}{2011}]{greenberg2011diagonally}
\begin{barticle}
\bauthor{\bsnm{Greenberg}, \binits{N.}},
\bauthor{\bsnm{Miller}, \binits{J.S.}}:
\batitle{Diagonally non-recursive functions and effective hausdorff dimension}.
\bjtitle{Bulletin of the London Mathematical Society}
\bvolume{43}(\bissue{4}),
\bfpage{636}--\blpage{654}
(\byear{2011})
\end{barticle}
\endbibitem

\bibitem[\protect\citeauthoryear{Liu}{2019}]{liu2019dnr}
\begin{botherref}
\oauthor{\bsnm{Liu}, \binits{L.}}:
Which dnr can be minimal.
arXiv preprint arXiv:1912.09053
(2019)
\end{botherref}
\endbibitem

\bibitem[\protect\citeauthoryear{Downey et~al.}{2002}]{downey2002trivial}
\begin{barticle}
\bauthor{\bsnm{Downey}, \binits{R.G.}},
\bauthor{\bsnm{Hirschfeldt}, \binits{D.R.}},
\bauthor{\bsnm{Nies}, \binits{A.}},
\bauthor{\bsnm{Stephan}, \binits{F.}}:
\batitle{Trivial reals}.
\bjtitle{Electronic Notes in Theoretical Computer Science}
\bvolume{66}(\bissue{1}),
\bfpage{36}--\blpage{52}
(\byear{2002})
\end{barticle}
\endbibitem

\bibitem[\protect\citeauthoryear{Nies}{2005}]{nies2005lowness}
\begin{barticle}
\bauthor{\bsnm{Nies}, \binits{A.}}:
\batitle{Lowness properties and randomness}.
\bjtitle{Advances in Mathematics}
\bvolume{197}(\bissue{1}),
\bfpage{274}--\blpage{305}
(\byear{2005})
\end{barticle}
\endbibitem

\end{thebibliography}

\end{document}